\documentclass[a4paper, 10pt]{amsproc}
\usepackage{multicol}
\usepackage{defs}
\usepackage{amsaddr}
\usepackage{mathtools}
\usepackage{amsmath,tikz}
\usetikzlibrary{matrix}
\usepackage{derivative}
\usepackage{algorithm} 
\usepackage{algorithmic}  
\usepackage[linesnumbered,ruled,lined,algo2e]{algorithm2e} 
\usepackage[numbered,framed]{matlab-prettifier}
\usetikzlibrary{arrows}
\usepackage{orcidlink}
\usepackage{caption}
\usepackage{graphicx}
\usepackage{array}
\usepackage{tabularray}
\usepackage{float}
\usepackage{booktabs}
\usepackage{array, multirow}
\pgfplotsset{compat=1.18}

\usepackage[font = {footnotesize, it}, labelfont = {bf, it}]{caption}

\hypersetup{final=true}

\title[\(\quito\)]{QuITO: Numerical software for constrained nonlinear optimal control problems}

\author[S. Ganguly]{\vspace{-6mm}Siddhartha Ganguly\,\orcidlink{0000-0003-2046-2061}} 
\address{\vspace{-5mm}
	\faGroup\ Systems \& Control Engineering\\
	\faUniversity\ IIT Bombay, Powai\\
 \faHome\ \url{https://sites.google.com/view/siddhartha-ganguly}\\
	\faMapMarker\ Mumbai 400076, India
}
\author[N. Randad, R. A. D'Silva]{\vspace{-6mm}Nakul Randad\,, Rihan Aaron D'Silva} 
\address{\vspace{-5mm}
	\faGroup\ Aerospace Engineering\\
	\faUniversity\ IIT Bombay, Powai\\
	\faMapMarker\ Mumbai 400076, India
}
\author[M. Raj]{\vspace{-6mm}Mukesh Raj S\,} 
\address{\vspace{-5mm}
	\faGroup\ Mechanical Engineering\\
	\faUniversity\ IIT Bombay, Powai\\
	\faMapMarker\ Mumbai 400076, India
}
\author[D. Chatterjee]{\vspace{-6mm}Debasish Chatterjee\,\orcidlink{0000-0002-1718-653X}} 
\address{\vspace{-5mm}
	\faGroup\ Systems \& Control Engineering\\
	\faUniversity\ IIT Bombay, Powai\\
 \faHome\ \url{http://www.sc.iitb.ac.in/~chatterjee}\\
	\faMapMarker\ Mumbai 400076, India
}
\thanks{%
	\faEnvelope\ \texttt{sganguly@iitb.ac.in}, \texttt{nakulrandad@gmail.com}, \texttt{dslvaaron24@gmail.com}, \texttt{mukeshrajs2711@gmail.com}, \texttt{dchatter@iitb.ac.in}
}
\thanks{Siddhartha Ganguly is supported by the PMRF grant RSPMRF0262, from the Ministry of Education, Govt.\ of India.}
\date{\DTMnow}

\begin{document}

\maketitle

\begin{abstract}
We introduce the MATLAB-based software \(\quito\) (Quasi-Interpolation based Trajectory Optimization) to numerically solve a wide class of constrained nonlinear optimal control problems (OCP). The solver is based on the \(\quito\) (the same abbreviation) algorithm, which is a \emph{direct multiple shooting} (DMS) technique that leverages a particular type of quasi-interpolation scheme for control trajectory parameterization. The software is equipped with several options for numerical integration, and optimization solvers along with a Graphical User Interface (GUI) to make the process of designing and solving the OCPs smooth and seamless for users with minimum coding experience. We demonstrate with two benchmark numerical examples the procedure to generate constrained state and control trajectories using \(\quito\).
\end{abstract}


\newcommand{\paramh}{h}
\newcommand{\paramD}{\mathcal{D}}

\section{Introduction}
\label{sec:intro}
Constrained dynamical optimization problems, more widely known as constrained optimal control problems (OCP) or trajectory optimization problems, arise in a variety of fields ranging from the engineering sciences and economics through the medical sciences. Optimal control theory provides several powerful tools for synthesizing constrained control trajectories optimizing certain objectives. Historically, the Pontryagin maximum principle (PMP) \cite{ref:liber} has been one of the primary driving engines behind the solution of constrained OCPs.
The chief bottlenecks of the PMP are the theoretical complications associated with the derivation of the necessary conditions for optimality, especially for a general class of constrained OCPs with state and/or mixed state-control/path constraints, and the associated numerical complexity thereof. 

\par In recent years, primarily due to the difficulties involved with numerical solutions of indirect methods, interest has shifted towards \emph{direct optimization} methods for computing numerically tractable and consistent optimal control trajectories. These direct methods transcribe the original (infinite-dimensional) OCP as a finite-dimensional nonlinear optimization problem post discretization. The methods are driven via the \emph{discretize-then-optimize} motto; naturally, they require an approximation scheme at the level of the control or state trajectories or both, from the ground up. Approximating the control and state trajectories by finitely parameterizing them in terms of appropriate basis or generating functions reduces the original infinite-dimensional optimization problem to a finite-dimensional nonlinear program (NLP). The ensuing NLP can be solved using off-the-shelf numerical solvers having suitable accuracy guarantees.
\par One of the earliest and well-established direct methods for solving constrained OCPs proceed via the \emph{direct shooting method} (DMS) \cite{ref:bock1984multiple}, which leverages control trajectory parameterization over the span of finitely many suitably chosen functions. Numerical algorithms that employ DMS using piecewise constant functions are MISER \cite{ref:software_MISER_3_update}, ACADO \cite{ref:software_acado}, etc. \textit{Direct collocation} techniques \cite{ref:ross-PSMreview-2012} constitute another popular class of direct methods employing an approximation scheme at the level of \emph{both} the state and the control trajectories. 
\emph{Pseudospectral/orthogonal collocation} methods are more recent additions to the class of direct collocation methods; they employ global/local orthogonal polynomials with Lagrange interpolation to approximate the states/control trajectories. Well-known software packages employing orthogonal collocation methods are ICLOCS2 \cite{ref:software_ICLOCS_II}, \(\mathbb{GPOPS}-\mathbb{I}\) \cite{ref:software_GPOPS_I}, and \(\mathbb{GPOPS}\) -- \(\mathbb{II}\) \cite{ref:software_GPOPS_II}. Table \eqref{tab:solver_spec} collects few of these software tools based on their underlying transcription method.

\begin{table}[htbp]
\caption{Transcription and a few available software}
\begin{tblr}{l|c}
	\hline[2pt]
		\SetRow{azure9}\SetCell[c=3]{c} Various optimal control software
& 3-2
& 3-3
 \\
	\hline[1pt]
		Transcription & Solver \\ 
 \hline
Multiple shooting & ACADO, MISER, DIRCOL \cite{ref:software_DIRCOL}\\
 \hline
Direct collocation & \(\mathbb{GPOPS}-\mathbb{I}\), \(\mathbb{GPOPS}-\mathbb{II}\), ICLOCS2  \\
	\hline[2pt]
\end{tblr}
\centering
\label{tab:solver_spec}
\end{table}


\subsection*{Contributions and impact overview}
\begin{itemize}[leftmargin=*, label=\(\circ\)]

	\item \(\quito\) is different from the other conventional solvers in the sense that the underlying approximation scheme is a \emph{quasi-interpolation} (instead of standard interpolation), and this enables the calculation of numerical control trajectories in a fast and matrix-free manner as opposed to other polynomial or radial basis function-based solvers and conventional collocation solvers \cite{ref:garg-2010, ref:software_ICLOCS_II} that requires matrix inversion. 

\item We record two benchmark numerical examples to demonstrate that for certain classes of optimal control problems \emph{our technique and the software performs better than some of the most advanced solvers currently available (which implement several advanced methods as specified earlier with possible mesh refinement modules)}. For an illustration, see the singular problems (singular optimal control problems are infamous \cite{ref:maurer2001sensitivity} for the numerical difficulty involved to solve them) in \eqref{numexmp:AlyChan} and \eqref{numexp:VanDerPol} where the control trajectory \emph{rings} \cite{ref:neuenhofen2018dynamic}, for both the cases, when the Legendre-Gauss-Radau (LGR) pseudospectral method was employed. \(\quito\), on the other hand, solves both problems seamlessly and produces solutions without any ringing. 

\item We also provide MATLAB code snippets for the numerical examples so that the users can easily implement their own OCPs using the \texttt{template.m} file in \(\quito\).

\end{itemize}
\(\quito\) supports the integration/discretization schemes and NLP solvers specified in Table \ref{tab:quito_specifications}. \(\quito\) employs CasADi \cite{ref:CASADI_andersson2019} to perform its symbolic declaration tasks and to interface with the interior point nonlinear optimization solver IPOPT \cite{ref:IPOPT_Biegler}, but other NLP and Quadratic Programming (QP) solvers compatible with CasADi can be readily included the package. 
      
           
         

\begin{table}[htbp]
\caption{Options available in \(\quito\).}
\begin{tblr}{l|c|c}
	\hline[2pt]
		\SetRow{azure9}\SetCell[c=4]{c} Discretization and solver options
& 3-2
& 3-3
& 3-4
\\
\hline[1pt]
Transcription & Discretization & NLP/QP \\ 
 \hline
\(\quito\) (DMS) & Euler, Trapezoidal, &  IPOPT, fmincon, and\\
& Runge-Kutta 4, Hermite Simpson & compatible QP solvers \\

	\hline[2pt]
\end{tblr}
\centering
\label{tab:quito_specifications}
\end{table}

\subsection*{Organization}

This article unfolds as follows. The primary problem that we address in this article is formulated in \S \eqref{sec:pre}. A comprehensive description of \(\quito\), its functionalities, architecture, and the graphical user interface (GUI) implementation details are given in \S\ref{sec:soft_dscr}. Two elaborate numerical examples are presented with an exhaustive comparison with the state-of-the-art software modules in \S\eqref{sec:numericalexp} to illustrate the effectiveness of the established algorithm and the software. Finally \S\ref{sec:conclu_dissc} outlines potential impacts and various industries and domains in engineering where \(\quito\) can be applied. 

\vspace{-2mm}

\subsection*{Notation}
We employ the standard notation. \(\N \Let \aset[]{1,2,\ldots}\) denotes the set of positive integers. We let \(\Nz\Let \Nzr\) denote the set of natural numbers and \(\Z\) denote the integers. The vector space \(\R[d]\) is assumed to be equipped with standard inner product \(\inprod{v}{v'}\Let \sum_{j=1}^d v_j v'_j\) for every \(v,v' \in \R[d]\).

\section{Preliminaries}\label{sec:pre}
We begin with formulating the class of constrained optimal control problems that \(\quito\) can solve. Let \(d,m \in \N\) and fix a final time horizon \(\horizon>0.\) Consider a nonlinear controlled dynamical system:
\begin{equation}
	\label{eq:sys}
	\dot \st(t) = \sys\bigl(\st(t)\bigr)+ G \bigl(\st(t)\bigr)\cont(t)\,\,\text{for a.e.\ } t \in \lcrc{0}{T}, 
\end{equation}
where \(\st(t) \in \R[d]\) is the vector of states at time \(t\), \(\cont(t) \in \Ubb \) is the control action at time $t$, where \(\Ubb\subset\Rbb^m\) is a given nonempty, compact set. The boundary states \(\bigl(\st(0), \st(\tfin)\bigr)\) are assumed to lie in a given nonempty closed subset of \(\tS \subset \Rbb^d\times\Rbb^d\). The control trajectory \(t\mapsto \cont(t)\) satisfies
\begin{equation}
	\label{eq:control constraint}
	\cont(\cdot)\in \mathcal{U}\Let \aset[\big]{u(\cdot) \in \lpL[\infty]\bigl(\lcrc{0}{\horizon};\mathbb{U})\suchthat u(t) \in \admcont \,\,\text{for a.e}\,\,t \in \lcrc{0}{T}}.
\end{equation}
A control \(\cont(\cdot)\) is \textit{feasible} if it is Lebesgue measurable, satisfies \eqref{eq:control constraint}, and the corresponding unique solution \(\st(\cdot)\) of \eqref{eq:sys} satisfies the boundary constraints. The acronym `a.e.' refers to almost everywhere relative to the Lebesgue measure; see \cite[Chapter 1]{ref:folland1999real} for technical details. We consider the OCP
\begin{equation}
\label{eq:OCP} 
\begin{aligned}
	& \minimize_{\cont(\cdot) \in \mathcal{U}} &&	\J(\st, \cont(\cdot)) \Let c_F(\st(\tfin))+ \int_{0}^{T}c(\st(t),\cont(t)) \, \dd t \\
&\sbjto &&  \begin{cases}
\text{dynamics}\,\, \eqref{eq:sys} \,\text{with its associated data},\\
\bigl(\st(0), \st(T)\bigr) \in \tS \subset\Rbb^d\times\Rbb^d,\\
	d_j\bigl(\st(t),\cont(t)\bigr) \le 0 \text{ for a.e}\,\,t \in \lcrc{0}{T},\\
\cont(t) \in \mathbb{U}\text{ for a.e }t \in \lcrc{0}{\horizon},
\end{cases}
\end{aligned}
\end{equation}
with the following data: 
\begin{itemize}[label=\(\circ\), leftmargin=*]
	\item The maps representing the drift vector field \(f: \Rbb^d \lra \Rbb^d\) and the control vector fields \(G: \Rbb^d \lra \Rbb^{d \times m}\) are continuous.
	\item The \emph{instantaneous cost} or the \emph{Lagrange cost}
		\(
			\Rbb^d \times \Rbb^m \ni (\dummyx,\dummyu) \mapsto c(\dummyx,\dummyu) \in \lcro{0}{+\infty}
		\)
		and the \emph{terminal cost} or the \emph{Mayer cost}
		\(
			\Rbb^d \ni \xi \mapsto c_F(\xi)\in \lcro{0}{+\infty}
		\)
		are continuous and strictly convex functions with compact sublevel sets.
	\item For each \(j=1,\ldots,R_0\), where \(R_0 \in \N\) the function \(\dummyx \mapsto d_j (\dummyx,\dummyu)\) is continuous and represents the mixed state-control constraints. 
\end{itemize}

\begin{remark}
	\(\quito\) was originally advanced under quadratic instantaneous cost functions (see \cite{ref:SG:NR:DC:RB-22}); however, empirical evidence suggests that \(\quito\) performs well under more general cost functions such as members of the aforementioned class. Convergence results of direct multiple shooting algorithms are typically derived under certain convexity and smoothness conditions on the problem data; see, e.g., \cite{ref:Bons-13} and the references therein. From a theoretical standpoint, a mathematical proof of uniform convergence of a multiple shooting algorithm for \emph{general classes} of optimal control problems, is open.
\end{remark}

\begin{remark}\label{rmk:lip_con}
Under certain smoothness assumptions on the cost functions and the dynamics, topological conditions on the constraint sets, and constraint qualification conditions it is possible to establish certain regularity results for the optimal control trajectory \(t \mapsto u\as(t)\). For example, see \cite{ref:depinho-2011lipschitz} and \cite{ref:vinter-galbraith-lipschitz} where Lipschitz continuity properties of \(u\as(\cdot)\) were established. 
\end{remark}

\section{Software description}\label{sec:soft_dscr}
We begin with a brief presentation of the technical engine driving our algorithm; more details on the algorithm may be found in \cite[\S III B]{ref:SG:NR:DC:RB-22} and mathematical background can be found in \cite[\S 7.2 and \S 7.3]{ref:GanCha-22}. We parameterize the control trajectory via the following one-dimensional approximation engine:
\begin{equation}\label{eq:app_app_on_R}
t \mapsto \overline{u}(t)\Let \frac{1}{\sqrt{\Dd }}\sum_{mh\in \fset}u(mh)\,\psi \biggl( \frac{t-mh}{h \sqrt{\Dd}}\biggr), \end{equation}
where \(h>0\) is the step size, \(\Dd>0\) is the shape parameter, \(\fset\) is a finite set, the generating function \(\psi(\cdot)\) is in the Schwartz class of functions \(\mathcal{S}(\Rbb;\Rbb)\), which need to satisfy a moment condition (of order \(N \in \N\)) and a decay condition, and \((u(mh))_{m\in \fset}\) are the unknown control coefficients to be determined by solving a finite-dimensional NLP. A host of choices \cite[Chapter 3]{ref:mazyabook} are available for the generating function \(\psi(\cdot)\); Table \ref{tab:basis_specifications} collects a few.
    
 




\begin{table}[htbp]
\caption{A list of generating functions and their corresponding order}
\begin{tblr}{l|c}
	\hline[2pt]
		\SetRow{azure9}\SetCell[c=3]{c} Generating functions \(\psi(\cdot)\) and order \(N \in \N\) of moment condition
& 3-2
& 3-3
 \\
	\hline[1pt]
		Generating functions & Order \\ 
 \hline
 \(\psi_1(x)=\frac{1}{\sqrt{\pi}}e^{-x^2}\)& \(2\)  \\ 
 \hline
\(\psi_2(x)=\frac{1}{\sqrt{\pi}}\left(\frac{3}{2}-x\right)e^{-x^2}\) & \(4\) \\
 \hline
\(\psi_3(x)=\frac{1}{\sqrt{\pi}}\left(\frac{15}{8}-\frac{5}{2}x^2+\frac{x^4}{4}\right)e^{-x^2}\) & \(6\) \\
 \hline
\(\psi_4(x)=\frac{1}{\sqrt{\pi}}\bigl(\frac{315}{128}-\frac{105}{16}x^2+\frac{63}{10}x^4-\frac{3}{4}x^6\) & \(10\) \\
\hspace{8mm}\(+\frac{1}{24}x^8\bigr)e^{-x^2}\) & \\
 \hline
\(\psi_5(x) = \sqrt{\frac{\mathrm{e}}{\pi}}e^{-x^2}\cos \sqrt{2}x\) & \(4\) \\ 
 \hline
\(\psi_6(x) = \frac{1}{\pi}\mathrm{sech} x\) & \(2\) \\
	\hline[2pt]
\end{tblr}
\centering
\label{tab:basis_specifications}
\end{table}

For the multiple shooting procedure, we solve the ODE \eqref{eq:sys} on each time interval \([mh,(m+1)h]\) with some initial guess and impose the \emph{continuity/defect condition} \cite[Eq. III.11]{ref:SG:NR:DC:RB-22} to ensure continuity of the trajectories. Finally, we arrive at an NLP with the discrete equivalent of all the constraints, and a quadrature approximation of the cost function \(\mathbb{J}(\cdot,\cdot)\), which can be solved using various off-the-shelf NLP solvers.
\subsection{Software architecture and functionalities}
\(\quito\) uses CasADi for symbolic variable declaration and interface with the optimization solvers, i.e., we call the NLP solver IPOPT and specify the constraints using CasADi's \texttt{opti} framework. For usage, we refer the readers to the official documentation at \url{https://web.casadi.org/docs/}. Other frameworks can also be employed instead of CasADi for modeling and optimization, for example, YALMIP \cite{ref:YALMIP_lofberg2004}. \(\quito\) consists of four core blocks:
\begin{enumerate}
    \item The directory \texttt{TemplateProblem} contains the template of the software following which the example problems are coded. Users can add their own problem easily by adding relevant problem data using the contents of the \texttt{TemplateProblem} folder. 

    \item The directory \texttt{examples} collects all the pre-coded examples that \(\quito\) consists of. The README.md files inside all the examples contain instructions on formulating the problems.

    \item The directory \texttt{src/problemTranscription} is one of the chief ingredients of \(\quito\) where the direct multiple shooting algorithm is encoded. From time specification, dynamics, adding constraints, solving the NLP, and generating approximate trajectories --- everything is encoded in file \texttt{solveProblem.m} within the \texttt{src/problemTranscription} folder. The directory \texttt{src/problemTranscription} must be added to the MATLAB path.

    \item The directory \texttt{Graphical Interface} builds the GUI for a one-click solution to the example problems; see \S\ref{subsec:GuI} and Figure \ref{fig:GUI_ss} for more details on the GUI implementation and a pictorial view of the GUI. 
\end{enumerate}
\(\quito\) builds the OCP by creating a problem template as given in the directory \texttt{TemplateProblem}, which must contain the files \texttt{TemplateProblem.m}, \texttt{main.m}, \texttt{options.m} and \texttt{postProcess.m}. The problem data needs to be defined in \texttt{TemplateProblem.m} by specifying the following: the continuous-time dynamics (set of ODEs), the state and the control dimensions, the Lagrange/instantaneous cost, the Mayer/terminal cost, the initial and final time, the state and control constraints, the set of initial conditions for states and the terminal state constraints. We set the default transcription technique as \texttt{options.transcription=QuITO} in  \texttt{options.m}. Preferred discretization methods \texttt{options.discretization}, choice of nonlinear program solver \texttt{options.NLPsolver} and corresponding solver parameters, generating functions \texttt{options.generating\_function}, meshing strategy (for now we employ uniform fixed meshing strategy) \texttt{options.meshstrategy}, the value of the shape parameter \(\Dd\) \texttt{options.variance}, computation time to solve the problem \texttt{options.print.time}, cost value \texttt{options.print.cost}, output plots \texttt{options.plot} etc. can all be modified in the \texttt{options.m} file as per the problem specification. 

Finally, to solve an OCP using \(\quito\) the user needs to choose two of the parameters related to the approximation scheme in the \texttt{main.m} file. In particular, the user needs to specify the arguments of the \texttt{solveProblem} command, i.e., the pair \((N,\Dd)\) where the first argument \(N\) is the total number of steps by which the time horizon \(\lcrc{0}{\horizon}\) has been divided, and based on a value of the step size \(h\) this number can be adjusted. The second argument is the shape parameter \(\Dd\), which can be modified in the command \texttt{options.variance} in the \texttt{options.m} file. The pair \((N,\Dd)\) are passed as arguments to the \texttt{options} function in \texttt{options.m}. The file \texttt{main.m} accesses the problem parameters and optimization specifications from files \texttt{problem.m} and \texttt{options.m} and computes the control trajectory and the state trajectory using the \texttt{solveProblem.m} file (available in the path) and then plots the same as specified in \texttt{postProcess.m}.
\begin{figure}[htbp]
\centerline{\includegraphics[width=10cm,height=7cm]{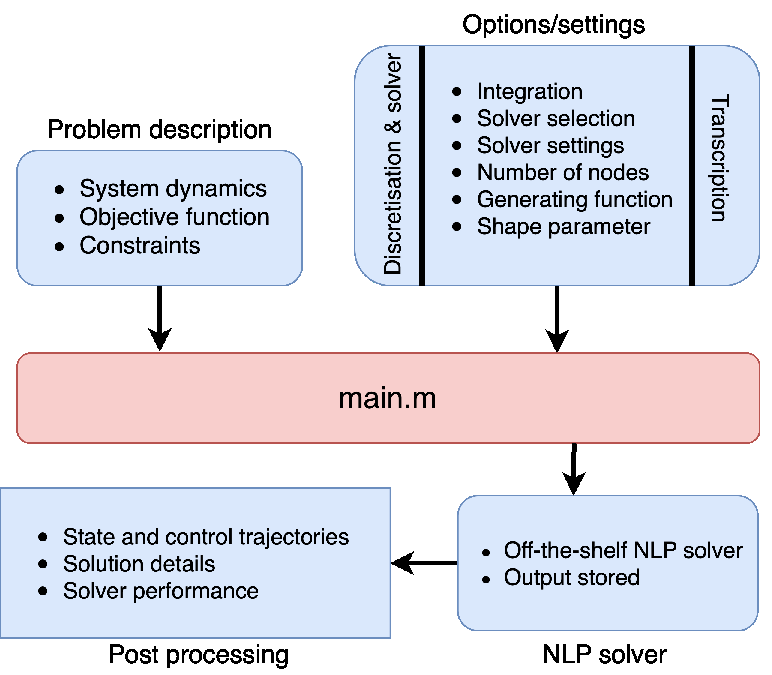}}
\caption{Flowchart of \(\quito\).}
\label{fig:QuITO_flowchart}
\end{figure}
In the unlikely case that the desired OCP does not fit into the default structure built into \(\quito\), the user can easily define the desired conditions (for example mixed state-control constraints) by including the same in a copy of a customized \texttt{solveProblem.m} file and place it in the problem directory. 

\subsubsection{General description}The flowchart in Figure \ref{fig:QuITO_flowchart} captures all essential components of the software and the interplay between them. Users can specify their problems in four different stages: specification of the problem data, objective and constraints falls under the \emph{Problem Description} block and these can be specified in \texttt{problem.m} file after that users can specify several types of discretization schemes, choice of basis functions, NLP solver setting in the \texttt{options.m} file, next the problem is transferred to the NLP solver which provides set of optimization control coefficients and finally the \texttt{postProcess.m} constructs the approximate control and state trajectories. All the blocks are linked to the \texttt{main.m} file which the user needs to run specifying the number of steps and the shape parameter of the generating functions.
\subsection{GUI implementation}\label{subsec:GuI}
\begin{figure}[htbp]
\centerline{\includegraphics[width=10cm,height=6cm]{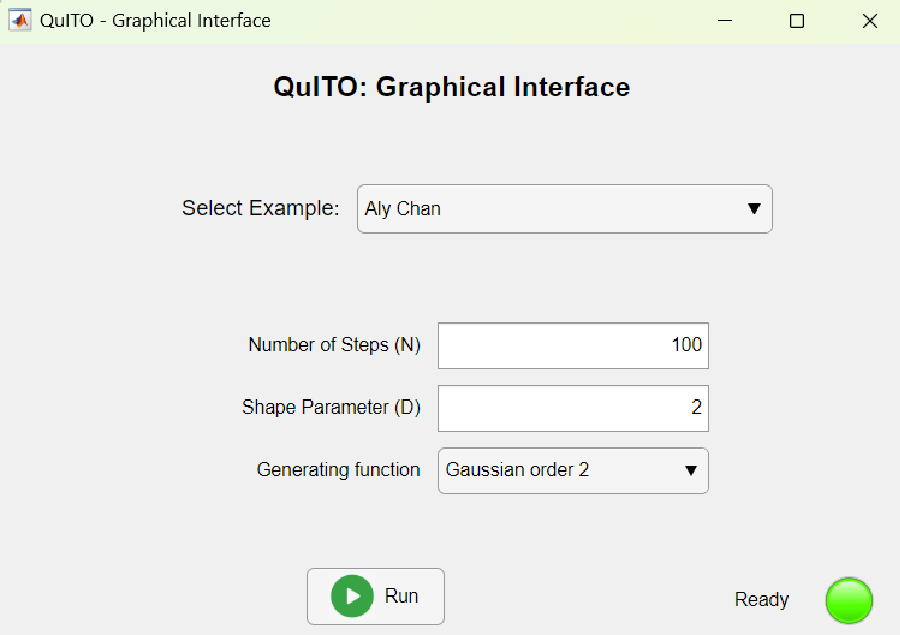}}
\caption{GUI screen for Example \ref{numexmp:AlyChan}.}
\label{fig:GUI_ss}
\end{figure}

To facilitate a more accessible and straightforward methodology to make use of the \(\quito\) toolbox, a linked MATLAB application with a Graphical User Interface (GUI) has been developed. Utilizing the GUI, users can choose any of the example OCPs currently available in \(\quito\) from a drop-down menu. The user can choose the parameters related to the approximation scheme \eqref{eq:app_app_on_R}, which include the number of steps, the generating function \(\psi(\cdot)\) and its shape parameter \(\Dd \in \loro{0}{+\infty}\); then simply click on the `Run' button to obtain the state and control trajectories associated with the chosen problem. Figure \ref{fig:GUI_ss} (right-hand sub-figure) shows the graphical interface for example \S\ref{numexmp:AlyChan}. 
The GUI interactively takes in inputs from the users for the parameters mentioned above and in turn calls the relevant MATLAB scripts that build and solve the chosen optimization problem and display the corresponding plots, \emph{all with the click of a single button}. The GUI implementation is contained in the Graphical Interface directory and in order to improve user accessibility, all requirements are packaged into a single file, \(\texttt{QuITO\_graphicalInterface.mlapp}\).  The front-end graphics of the GUI was developed making use of MATLAB's App Designer. The back-end processing was built to work in the following sequence: 
\begin{enumerate}
    \item The user picks the \(\quito\) parameters for the problem they wish to solve and clicks on the `Run' button, consequently  MATLAB navigates to the relevant directory and appends it to MATLAB's searchable path. 
    \item Problem description for the chosen problem, describing the system dynamics, objective function, and constraints, is loaded and the options file containing the appropriate options to solve the problem, as chosen by the user, is generated.
    \item Utilising the preceding files the solver is initiated to solve the problem. Once the problem is solved, the solution is loaded onto a structure after which the solution is post-processed and appropriate plots are displayed.
    \item Then, the changes made earlier to MATLAB's path are cleared and the system is prepared for further operations and waits for the next set of inputs from the user. 
\end{enumerate}

The GUI additionally enables solving any template problem the user desires. The user would define their problem by populating the files in the \texttt{TemplateProblem} directory given in the package. Then, open the MATLAB application \(\texttt{QuITO\_graphicalInterface.mlapp}\), select \texttt{User-defined template problem} from the drop-down, choose the parameters \((\psi(\cdot),N,\Dd)\) and click \texttt{Run} to obtain the results. Detailed instructions for running the GUI can be found on the package's GitHub repository. 
\section{Numerical experiments} \label{sec:numericalexp}
This section contains two numerical experiments showcasing the superior performance of \(\quito\) on benchmark systems over and above the current state-of-the-art. Both examples demonstrate that in the absence of adaptive mesh refinement, \(\quito\) performs better than state-of-the-art algorithms for benchmark optimal control problems. We compare our results with the state-of-the-art collocation techniques via the ICLOCS2 solver. The software package can be found on the GitHub repository: 

\begin{table}[htbp]

\begin{tblr}{l}
	\hline[2pt]
		\SetRow{red9}\SetCell[c=3]{c} \url{https://github.com/Gdarthsid/QuITO-Version-1.git}
& 3-2
& 3-3
 \\
	\hline[2pt]

\end{tblr}
\centering
\label{solver_spec}
\end{table}
\subsection{Aly-Chan problem} \label{numexmp:AlyChan}
\begin{figure}[htbp]
\centerline{\includegraphics[width=10cm,height=8cm]{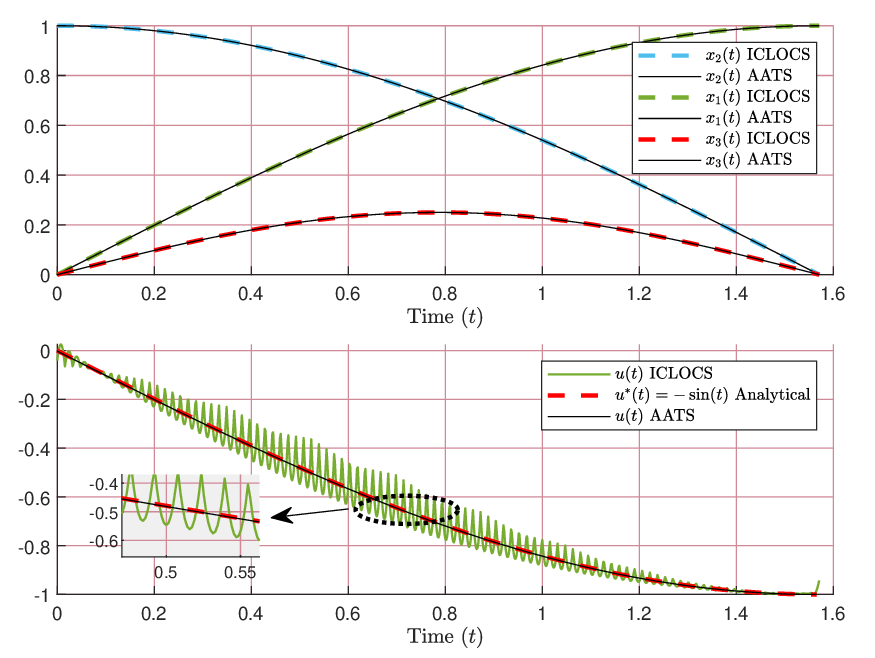}}
\caption{State and control trajectories for Example \ref{numexmp:AlyChan}.}
\label{fig:traj_AC}
\end{figure}

\begin{figure}[ht]
\includegraphics[width=0.48\textwidth]{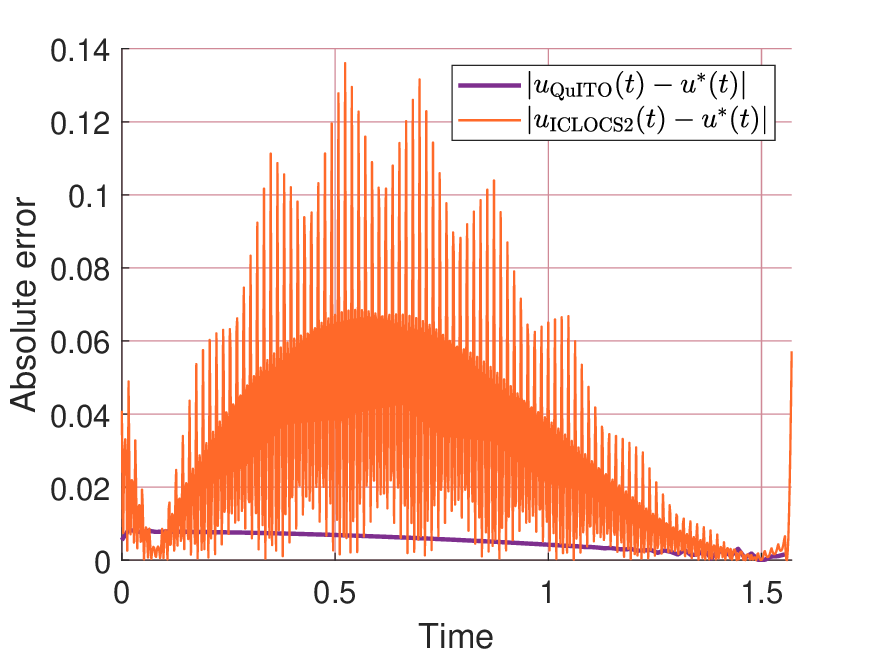}
\hspace{\fill}
\includegraphics[width=0.48\textwidth]{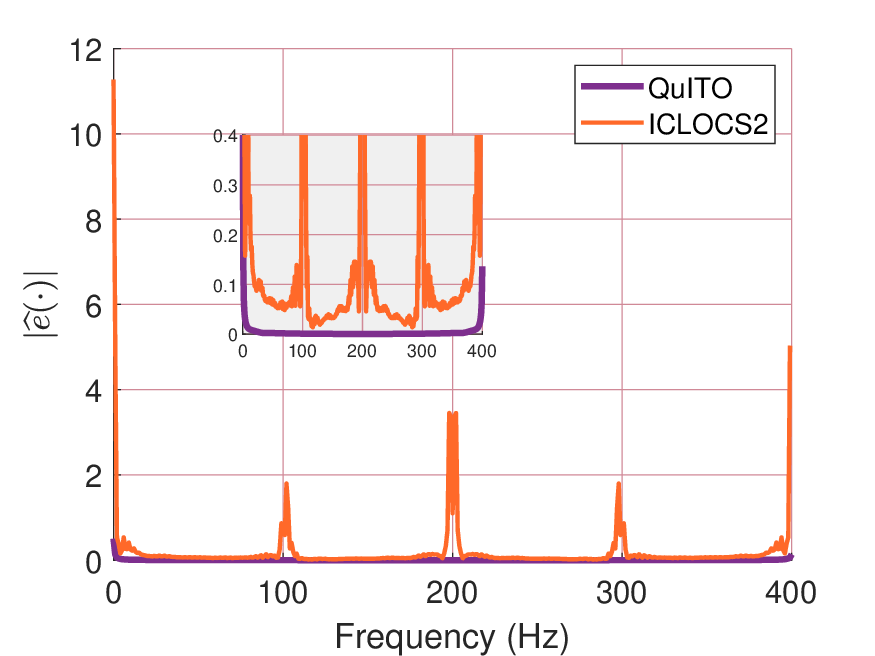}
\caption{Control error trajectories for Example \ref{numexmp:AlyChan} and their corresponding DFT plots with \(\quito\) and ICLOCS2 respectively.}\label{fig:AC_error_dft}
\end{figure}

Consider the so-called \emph{Aly-Chan} OCP with third-order dynamics and hard constraints on control:
\begin{equation}
\label{eq:aly-chan-ocp}
\begin{aligned}
& \minimize_{\cont(\cdot)}	&&  x_3 \left( \tfin\right) \\
&  \sbjto		&&  \begin{cases}
\dot{\st}_1(t) =x_2(t),\, \dot{\st}_2(t)=u(t),\\  \dot{\st}_3(t)=\frac{1}{2}\left(x_2(t)^2-x_1(t)^2\right),\\ 
\st_1(0)=0, \,\st_2(0)=1, \,\,\text{and}\,\,\st_3(0)=0,\\
|u(t)|\le 1, \,\, \tfin = \tfrac{\pi}{2}.
\end{cases}
\end{aligned}
\end{equation}
The analytical expression of the optimal control trajectory is known and is given by \(t \mapsto u\as(t) = -\sin(t)\). The optimal control trajectory exhibits a singular behaviour; we refer to \cite{ref:aly-chan} for a detailed treatment. We provide a brief overview of the design procedure of the Aly-Chan problem in \(\quito\). We start with defining the system dynamics and subsequently specify the objective function and constraints respectively. All of these are encoded in \texttt{AlyChan.m}:
\begin{lstlisting}[
  style      = Matlab-editor,
  basicstyle = \mlttfamily,
  numbers = none,
]
function [dx] = dynamics(x,u,t)
dx1 = x(2); dx2 = u(1); dx3 = 0.5*(x(2)^2 - x(1)^2);
dx = [dx1; dx2; dx3]; end
\end{lstlisting}
Next, we encode the objective function, which consists of a Mayer/terminal cost:  
\begin{lstlisting}[
  style      = Matlab-editor,
  basicstyle = \mlttfamily,
   numbers = none,
]
function lag = stageCost(x,u,t) 
lag = 0; end
function mayer = terminalCost(x,u,t)
mayer = x(3); end
\end{lstlisting}
Next, we call the functions \texttt{dynamics}, \texttt{stageCost}, and \texttt{terminalCost}:
\begin{lstlisting}[
  style      = Matlab-editor,
  basicstyle = \mlttfamily,
   numbers = none,
]
problem.dynamicsFunc = @dynamics; 
problem.stageCost = @stageCost; 
problem.terminalCost = @terminalCost;
\end{lstlisting}
and specify the system parameters and the constraints: 
\begin{lstlisting}[
  style      = Matlab-editor,
  basicstyle = \mlttfamily,
   numbers = none,
]
problem.time.t0 = 0; % initial time
problem.time.tf = pi/2; % final time
problem.nx = 3; % Number of states
problem.nu = 1; % Number of controls
problem.states.x0l = [0 1 0]; 
problem.states.x0u = [0 1 0]; 
problem.states.xl = [-inf -inf -inf];
problem.states.xu = [inf inf inf];
problem.states.xfl = [-inf -inf -inf]; 
problem.states.xfu = [inf inf inf];
problem.inputs.ul = -1;
problem.inputs.uu = 1;
\end{lstlisting}
By equating \texttt{options.generating\_function} to \texttt{generating\_function\_flag} 
we set the generating function for the quasi-interpolation in \texttt{options.m} in directory \texttt{Aly Chan}. By default, \texttt{generating\_function\_flag} takes \(\psi_{1}(x) \Let e^{-x^2}/\sqrt{\pi}\) as the generating function. The user can choose other options by setting \texttt{generating\_function\_flag} to a different value between \(2\) and \(6\). For discretization and NLP solve we picked the Euler discretization and the solver IPOPT respectively: \texttt{options.discretization=`euler'} and \texttt{options.NLPsolver = `ipopt'} in \texttt{options.m} inside the \texttt{Aly Chan} folder, which is called along with \texttt{problem.m} in the \texttt{main.m} file. In \texttt{main.m} we specify the pair \((N,\Dd)=(100,2)\). To plot the state and control trajectories we choose \texttt{options.plot = 2} in \texttt{options.m}: 
\begin{lstlisting}[
  style      = Matlab-editor,
  basicstyle = \mlttfamily,
   numbers = none,
]
problem = AlyChan; 
opts = options(100, 2);
solution = solveProblem(problem, opts);
postProcess(solution, problem, opts)
\end{lstlisting}
It was observed that the optimal control trajectory displayed an unwanted behaviour---the so-called \emph{ringing phenomenon} \cite{ref:neuenhofen2018dynamic} when the LGR pseudospectral method was employed via the package ICLOCS2. In contrast, the \(\quito\) solution exhibits no such peculiarities. Figure \ref{fig:traj_AC} shows the analytical solution, the solution obtained from ICLOCS2 (via LGR collocation) and \(\quito\), all on the same time horizon  \(\lcrc{0}{\frac{\pi}{2}}\), while Figure \ref{fig:GUI_ss} shows the GUI interface for the Aly Chan problem. The control error trajectories and their corresponding DFT (Discrete Fourier Transform) plots are given in Figure \ref{fig:AC_error_dft} and Figure \ref{fig:log_error} (left-hand sub-figure, in \emph{log scale}) when \(\quito\) and ICLOCS2 were employed, respectively. The \(\lpL[2]\)-norm \(\norm{e(\cdot)}_{\lpL[2]}\) of the error trajectories obtained from \(\quito\) and ICLOCS2 were, respectively, \(0.02\) and \(0.7816\) numerically computed via the command \texttt{norm(}\(\cdot\)\texttt{,2)}; moreover, the \(\ell^2\)-norm \(\norm{\widehat{e}(\cdot)}_{\ell^2}\) of the DFTs of the corresponding error trajectories were found to be \(0.558\) and \(15.632\), respectively, via standard MATLAB command \texttt{fft}\footnote{See the url \url{https://in.mathworks.com/help/matlab/ref/fft.html}.} and \texttt{norm}. Observed through the lens of time-frequency analysis, the preceding figures demonstrate that significantly less energy --- one order of magnitude --- is present in both the time- and frequency-spaces of the error trajectories (relative to the analytical solutions) due to \(\quito\) compared to those due to ICLOCS2. Table \ref{AC:table-cpu-time} collects certain optimization statistics for Example \ref{numexmp:AlyChan}.

\begin{table}[htbp]
\caption{A comparison of optimization statistics for Example \ref{numexmp:AlyChan}.}
\begin{tblr}{l|c|c|c}
	\hline[2pt]
		\SetRow{azure9}\SetCell[c=4]{c} CPU time, number of iterations and constraint violation comparison
& 3-2
& 3-3
& 3-4
\\
\hline[1pt]
Method & CPU time & Iteration count & Constraint violation \\ 
\hline
\(\quito\) & 4 sec  & 19   & \(1.945 \times 10^{-16}\)\\
 \hline
ICLOCS2 (optimal control rings) & 2 sec  & 18   & \(5.339 \times 10^{-10}\)\\
 \hline[2pt]
\end{tblr}
\centering
\label{AC:table-cpu-time}
\end{table}

\subsection{Van der Pol Oscillator}\label{numexp:VanDerPol}
We consider the Van der Pol oscillator system with bounded control:
    \begin{equation}
\label{eq:VanDerPol_OCP_con}
\begin{aligned}
& \minimize_{\cont(\cdot)}	&&  \frac{1}{2}\int_{0}^{4} \left(x_1(t)^2+x_2(t)^2\right)\, \dd t \\
&  \sbjto		&&  \begin{cases}
\dot{\st}_1(t)=\st_2(t),\\
\dot{\st}_2(t)=\bigl(1-\st_1^2(t)\bigr)\st_2(t)-\st_1(t)+\cont(t),\\
\bigl(\st_1(0)\,\,\st_2(0)\bigr)^{\top}=(0\,\,1)^{\top},\,\,|u(t)| \le 1.
\end{cases}
\end{aligned}
\end{equation}
\begin{figure}[tbp]
\centerline{\includegraphics[scale=0.7]{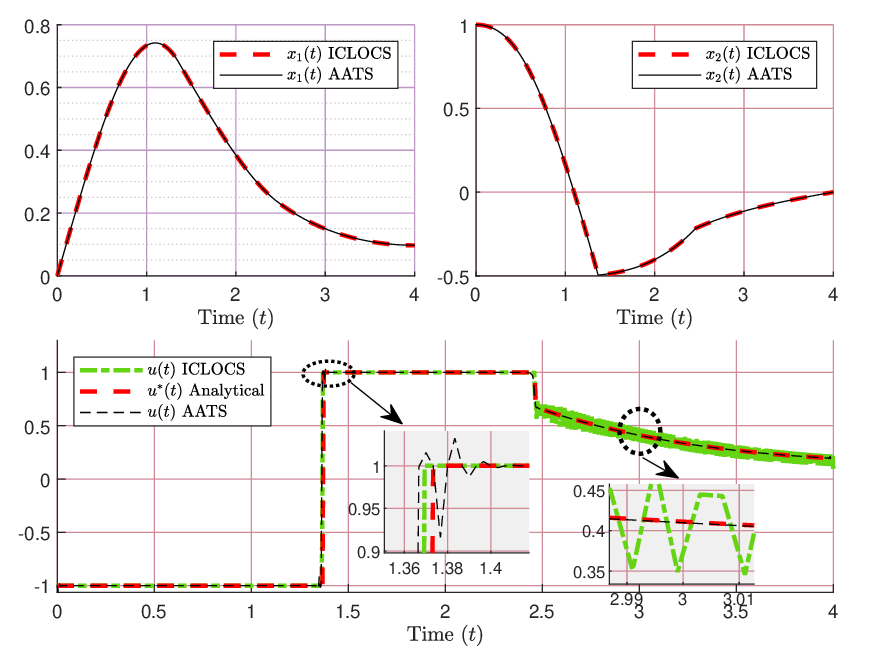}}
\caption{State and control trajectories for Example \ref{numexp:VanDerPol}.}
\label{fig:VanPol_st_con}
\end{figure}

\begin{figure}[ht]
\includegraphics[width=0.48\textwidth]{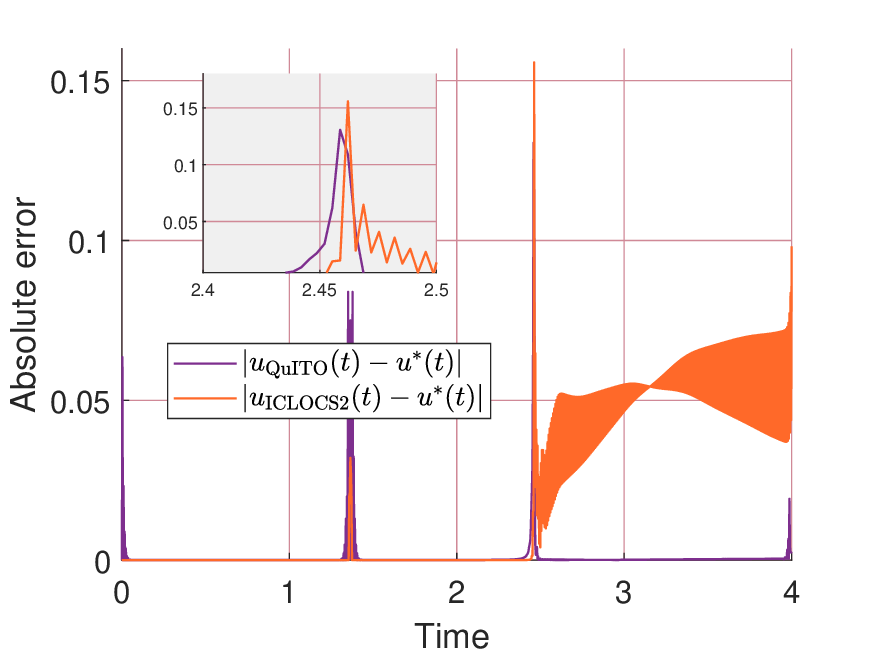}
\hspace{\fill}
\includegraphics[width=0.48\textwidth]{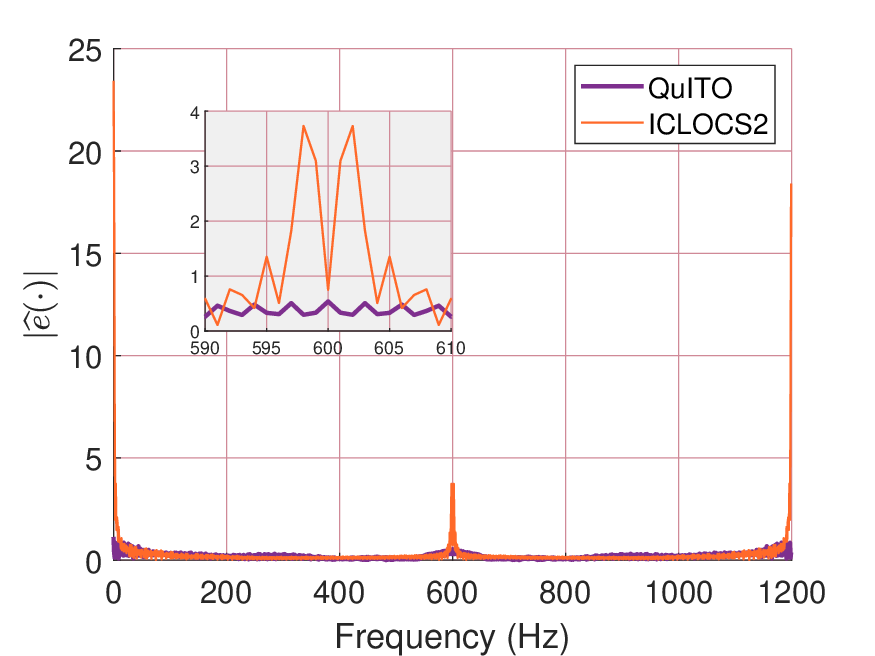}
\caption{Control error trajectories for Example \ref{numexp:VanDerPol} and their corresponding DFT plots with \(\quito\) and ICLOCS2 respectively.}\label{fig:VD_error_dft}
\end{figure}

\begin{figure}[ht]
\includegraphics[width=0.48\textwidth]{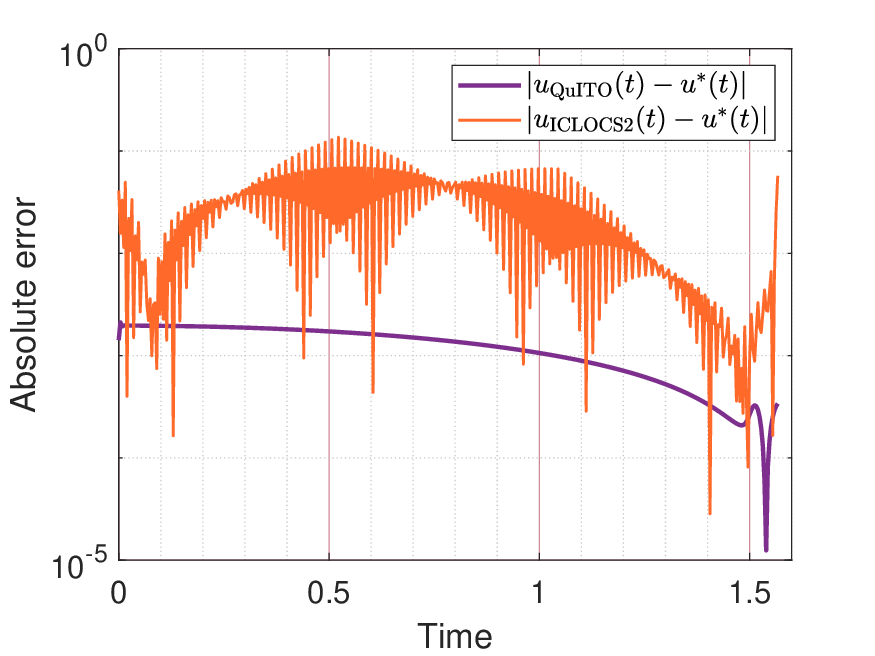}
\hspace{\fill}
\includegraphics[width=0.48\textwidth]{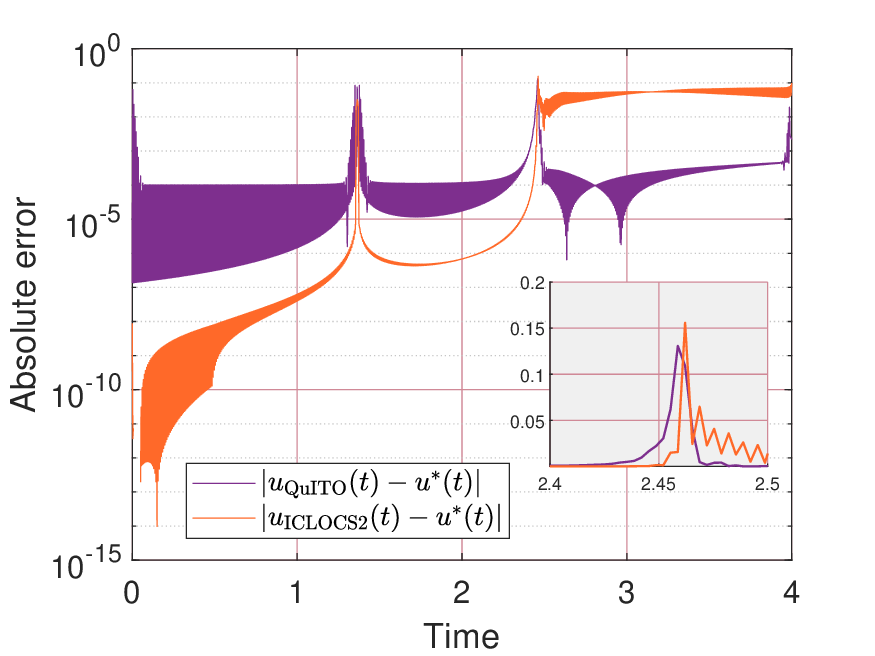}
\caption{Control error trajectories for Example \ref{numexmp:AlyChan} (left-hand sub-figure) and Example \ref{numexp:VanDerPol} (right-hand sub-figure) in \emph{log scale} with \(\quito\) and ICLOCS2 respectively.}\label{fig:log_error}
\end{figure}

Brief code snippets are given below: 
\begin{lstlisting}[
  style      = Matlab-editor,
  basicstyle = \mlttfamily,
  numbers = none,
]
function [dx] = dynamics(x,u,t)
dx1 = x(2); dx2 = (1 - x(1)^2)*x(2) - x(1) + u(1);
dx = [dx1; dx2]; end
\end{lstlisting}
The OCP \eqref{numexp:VanDerPol} consists of a Lagrange cost: 
\begin{lstlisting}[
  style      = Matlab-editor,
  basicstyle = \mlttfamily,
   numbers = none,
]
function lag = stageCost(x,u,t) 
lag = 0.5*(x(1)^2 + x(2)^2); end
function mayer = terminalCost(x,u,t)
mayer = 0; end
\end{lstlisting}
Next, we call the functions \texttt{dynamics}, \texttt{stageCost}, and \texttt{terminalCost}:
\begin{lstlisting}[
  style      = Matlab-editor,
  basicstyle = \mlttfamily,
   numbers = none,
]
problem.dynamicsFunc = @dynamics; 
problem.stageCost = @stageCost; 
problem.terminalCost = @terminalCost;
\end{lstlisting}
and specify the system parameters and the constraints: 
\begin{lstlisting}[
  style      = Matlab-editor,
  basicstyle = \mlttfamily,
   numbers = none,
]
problem.time.t0 = 0; % initial time
problem.time.tf = 4; % final time
problem.nx = 2; % Number of states
problem.nu = 1; % Number of controls
problem.states.x0l = [0 1]; 
problem.states.x0u = [0 1]; 
problem.states.xl = [-inf -inf];
problem.states.xu = [inf inf];
problem.states.xfl = [-inf -inf]; 
problem.states.xfu = [inf inf];
problem.inputs.ul = [-1];
problem.inputs.uu = [1];
\end{lstlisting}
Finally we execute \texttt{main.m} script: 
\begin{lstlisting}[
  style      = Matlab-editor,
  basicstyle = \mlttfamily,
   numbers = none,
]
problem = AlyChan; 
opts = options(300, 2);
solution = solveProblem(problem, opts);
postProcess(solution, problem, opts)
\end{lstlisting}
to obtain the state and control trajectories. The optimal control for this problem is known to feature discontinuities with a singular arc on a sub-interval of \([0,4]\) \cite{ref:neuenhofen2018dynamic}. The \emph{ringing phenomenon} (as noticed in the Aly-Chan problem in \S\ref{numexmp:AlyChan}; see the left-hand sub-figure of Figure \ref{fig:GUI_ss}) can be observed when the LGR collocation method was employed via ICLOCS2 with 400 nodes, while the \(\quito\) algorithm numerically provided better results; see Figure \ref{fig:VanPol_st_con}, Figure \ref{fig:VD_error_dft}, and Figure \ref{fig:log_error}. The \(\lpL[2]\)-norm \(\norm{e(\cdot)}_{\lpL[2]}\) of the error trajectories obtained from \(\quito\) and ICLOCS2 were numerically computed via the command \texttt{norm(}\(\cdot\)\texttt{,2)} to be, respectively \(0.2522\) and \(1.1327\); moreover, the \(\ell^2\)-norm \(\norm{\widehat{e}(\cdot)}_{\ell^2}\) of the DFT of the corresponding error trajectories were numerically computed via the standard MATLAB command \texttt{fft} and \texttt{norm} to be \(8.73\) and \(39.2371\). Observed through the lens of time-frequency analysis, the preceding figures establish that less energy is present in both the time- and frequency-spaces of the error trajectories (relative to the analytical solutions) obtained from \(\quito\) compared to those obtained from ICLOCS2.


	The preceding observations for two benchmark examples (Example \ref{numexmp:AlyChan} and Example \ref{numexp:VanDerPol}) indicate better quality of performance of \(\quito\) relative to ICLOCS2. For singular control problems, one often makes a compromise between the speed of computation and the accuracy of the solution. In this context, \emph{our main focus is the accuracy of the solution trajectories which is arguably the most important metric for accuracy- and safety-critical applications}. Table \ref{vdo:table-cpu-time} collects certain optimization statistics for Example \ref{numexp:VanDerPol}. Notice that compared to Example \ref{numexmp:AlyChan}, the computation time taken by \(\quito\) to generate its approximate control trajectory for Example \ref{numexp:VanDerPol} is higher than that by ICLOCS2 due to a fine discretization employed by \(\quito\) to capture the discontinuity. On the one hand, a fine discretization is necessary in \(\quito\) at this stage in the absence of mesh-refinement. On the other hand, a fine discretization is not required in LGR collocation methods because they employ piecewise constant generating functions to parameterize the control trajectory, and such a parametrization reduces the computation time. While these techniques pay attention to reducing the speed of the computation, \(\quito\) focuses on solution accuracy (without compromising too much on computation speed).

\begin{table}[htbp]
\caption{A comparison of optimization statistics for Example \ref{numexp:VanDerPol}.}
\begin{tblr}{l|c|c|c}
	\hline[2pt]
		\SetRow{azure9}\SetCell[c=4]{c} CPU time, number of iterations and constraint violation comparison
& 3-2
& 3-3
& 3-4
\\
\hline[1pt]
Method & CPU time & Iteration count & Constraint violation \\ 
\hline
\(\quito\) & 40 sec  & 14   & \(6.969 \times 10^{-14}\)\\
 \hline
 ICLOCS2 with mesh-refinement \\(optimal control rings) & 5 sec  & 16   & \(5.462 \times 10^{-19}\)\\
 \hline[2pt]
\end{tblr}
\centering
\label{vdo:table-cpu-time}
\end{table}
The (qualitative) oscillatory behaviour near the jumps can be reduced to some extent by employing selective local mesh-refinement in \(\quito\); theoretical results and their consequent software updates in this direction are currently under development and will be reported subsequently.

\section{Impact and Conclusion}\label{sec:conclu_dissc}
We introduced the software \emph{\(\quito\) --- Quasi-Interpolation based Trajectory Optimization} for numerically solving a wide class of nonlinear constrained optimal control and trajectory optimization problems \cite{ref:QuITO:SoftX}. It is based on the direct multiple shooting algorithm \(\quito\). We illustrated with two benchmark examples the solid performance of \(\quito\) despite it being currently at a nascent stage; indeed, \(\quito\) performed better than conventional solvers on important benchmark problems as demonstrated above. We specifically observed that for several classes of problems where state-of-the-art pseudospectral collocation generates control trajectories that \emph{ring}, \(\quito\) shows no such undesirable artifacts. The preceding data point to its remarkable potential for applications to a range of singular problems.

In \(\quito\), the user is only required to specify their problem data, choice of discretization, solver, and settings, plotting requirements in the template file \texttt{TemplateProblem.m}, \texttt{options.m}, and \texttt{postProcess.m} in the directory \texttt{TemplateProblem} and run \texttt{main.m} or alternatively run the GUI specific to the problem. The GUI interface of \(\quito\) is an especially attractive feature, which enables users to solve a specific problem quickly in a hassle-free manner, and makes it appealing for industry applications.

\(\quito\) is not restricted to optimal control problems; it can be used to solve a large class of practical optimization problems under constraints as long as the problem can be reformulated in the form \eqref{eq:OCP}. For example, various industries such as chemical process, automotive, robotics, aerospace, etc., essentially formulate their practical problems in the language of dynamic optimization, and consequently, \(\quito\) readily fits such applications. \(\quito\) contains a few benchmark industrial control systems, for example, an inverted pendulum on a cart, robotic motion planning and obstacle avoidance, etc. In fact, \(\quito\) \emph{has been employed in practice} \cite{ref:QuITO:patent} for a motion planning and obstacle avoidance problem for a robotics application and has performed significantly better than other competing methods available today; see the example \texttt{Robot Path Planning} in the \texttt{examples} directory for a numerical study. We are currently in the process of making the software package \(\quito\) computationally faster and more efficient by incorporating additional functionalities such as more solver options, discretization schemes, mesh refinement algorithms, etc., and subsequently including them to the GUI. Our future research direction includes making \(\quito\) more optimized and faster along with implementing local mesh-refinement for solving problems where the control trajectory is nonsmooth and even potentially discontinuous.
\section{Acknowledgement}
SG is supported by the PMRF grant RSPMRF0262, from the MHRD, Govt.\ of India. The authors thank Ravi Banavar from the Systems and Control Engineering Department, IIT Bombay, India for his helpful comments on the initial version of this manuscript and Pranav Lad from MathWorks India for his guidance during the initial phase of the GUI implementation.   
\bibliographystyle{amsalpha}
\bibliography{refs}

\end{document}